\documentclass[12pt]{article}
\usepackage[english]{babel}
\usepackage{amsmath,amsthm}
\usepackage{amsfonts,amssymb}
\usepackage[T1]{fontenc}
\usepackage[utf8]{inputenc}
\usepackage{indentfirst}
\usepackage{makeidx}
\usepackage[nottoc]{tocbibind}
\usepackage{courier}
\usepackage{type1cm}
\usepackage{listings}
\usepackage{titletoc}
\usepackage[all,cmtip]{xy}
\usepackage[pdftex]{graphicx}

\def\bes{\begin{eqnarray*}}
\def\ees{\end{eqnarray*}}
\def\bee{\begin{eqnarray}}
\def\eee{\end{eqnarray}}
\def\la{\langle}
\def\ra{\rangle}
\def\N{\mathbf N}

\def\C{\mathbf C}
\def\R{\mathbf R}

\def\a{\alpha}
\def\b{\beta}
\def\g{\gamma}

\def\s{\sigma}
\def\f{\varphi}

\def\Z{\mathbf Z}
\def\0{\bar 0}
\def\1{\bar 1}
\def\V{{\mathcal V}}

\def\prf{{\it Proof.\ } }
\def\ctd{\hfill$\Box$}

\def\VV{\mathcal V}

\newtheorem{thm}{Theorem}[section]

\newtheorem{lem}[thm]{Lemma}

\newtheorem{conj}[thm]{Conjecture}

\theoremstyle{definition}

\numberwithin{equation}{section}
\title{Simple Jordan superalgebras with the even parts of Clifford type}

\author{Ivan Shestakov\thanks{Supported  by FAPESP, grant 2018/23690-6,  by CNPq grant 304313/2019-0, and by IMC SUSTech, Shenzhen, China} \\
{\small Instituto de Matem\'atica e Estat\'{\i}stica}\\
{\small Universidade de S\~ao Paulo}\\
{\small S\~ao Paulo. Brasil}\\
{\small ivan.shestakov@gmail.com}
\and
Efim Zelmanov\thanks{The author gratefully acknowledges support from the NSF of China grant {No\,12350710787}}\\
{\small Shenzhen International Center for Mathematics,}\\
 {\small Southern University of Science and Technology, }\\
 {\small Shenzhen, China}\\
 {\small efim.zelmanov@gmail.com}
}

\date{\quad}

\begin{document}
\maketitle

\begin{abstract}
The purpose of this paper is a  partial progress towards classification of simple infinite dimensional Jordan superalgebras.   First, we prove that the only simple infinite dimensional Jordan superalgebras with  finite dimensional even parts are the superalgebras of superforms.
Then we consider the superalgebras whose even parts are infinite dimensional algebras of ``Clifford type'',  that is,   direct sums of algebras of bilinear forms.   The results of \cite{RZ} show that the number of summonds in these sums is 1 or 2.   We prove that the second case is impossible and that the  simple infinite dimensional Jordan superalgebras of the first type are the superalgebras of superforms.
\end{abstract}

% ----------------------------------------------------------------

\section{Definitions, examples, and conjectures}

\hspace{\parindent}
Throughout the paper all algebras are considered over a field $F$ of characteristic $\neq 2$.

A (linear ) {\em Jordan algebra} is a vector space $J$ with a binary bilinear operation $xy$ satisfying the identities
\bee
  xy&=&yx,\\
(x^2y)x&=&x^2(yx).
\eee

Let $V$ be a vector space of countable dimension and let $G(V)$ be the Grassmann  (or exterior) algebra over $V$. If $\{e_i,\ i\geq 1\}$ is a basis of the space $V$, then the algebra $G(V)$ is presented by generators $e_i,\ i\geq 1$,  and relations $e_ie_j+e_je_i=0,\ i.j\geq 1$. The set of ordered products $1,\, e_i,\, \ldots, e_{i_1}\cdots e_{i_k},\ i_1<i_2<\cdots<i_k,$  is a basis of $G(V)$.

The algebra $G(V)$ is $\Z/2\Z$-graded,  $G(V)=G(V)_{\0}+G(V)_{\1}$. The even part $G(V)_{\0}$ is the span of products of even length: $1,e_{i_1}\cdots e_{i_{2k}},\ i_1<\cdots <i_{2k}$, whereas $G(V)_{\1}$ is the space of all products of odd length $e_{i_1}\cdots e_{i_{2k+1}},\ i_1<\cdots<i_{2k+1}$.

By a {\em superalgebra} we mean an algebra
\bes
A=A_{\0}+A_{\1}
\ees
that is $\Z/2\Z$-graded.   The elements from $A_{\0}$ are usually called {\em even} and the elements from $A_{\1}$ are called {\em odd}.

Given a variety $\VV$ of algebras defined by homogeneous identities (see \cite{Jac, ZSSS}) we say that a superalgebra $A=A_{\0}+A_{\1}$ is a $\VV$-superalgebra if its {\em Grassmann envelope} 
\bes
G(A)=A_{\0}\otimes G(V)_{\0}+A_{\1}\otimes G(V)_{\1}
\ees
belongs to the variety $\VV$.

Knowing the identities defining the variety $\VV$,  one can easily write down the {\em superidentities}  defining the $\VV$-superalgebras.
Thus an associative superalgebra is just a $\Z/2\Z$-graded associative algebra.  A {\em commutative superalgebra}  is a $\Z/2\Z$-graded algebra $A=A_{\0}+A_{\1}$ that satisfies 
the superidentity 
\bes
xy=(-1)^{|x||y|}yx,
\ees
where for an element $a\in A_{\0}\cup A_{\1}$ the symbol  $|a|$ denotes its parity, $|a|=0$ or 1.
In particular,  the Grassmann algebra is a commutative superalgebra.
\smallskip

Observe that the field $\C=\R+\R i$  of compex numbers is a real superalgebra which is commutative as an algebra but not commutative as a  superalgebra.

Given an algebra $A$, one can construct the {\em $A$-double superalgebra}  $A(\sqrt 1)=A+Au$,  where $u(=\sqrt 1)$ is an odd central element such that $u^2=1$. We have $(A(\sqrt 1))_{\0}=A,\ (A(\sqrt 1))_{\1}=Au$.  One can easily prove that any simple superalgebra is either simple as an algebra or is isomorphic to $A(\sqrt 1)$ for some simple algerbra $A$. In particular, every simple  commutative superalgebra $A$  is simple as an algebra (since $x^2=0$ for any element $x\in A_{\1}$).
\smallskip

A superalgebra $J=J_{\0}+J_{\1}$ is a {\em Jordan superalgebra} if it satisfies the  superidentities
\bee
& xy=(-1)^{|x||y|}yx,&\label{SJ1}\\
& ((xy)z)t+(-1)^{|y||z|+|y||t|+|z||t|}((xt)z)y+(-1)^{|z||t|}x((yt)z)&\nonumber\\
&= (xy)(zt)+(-1)^{|y||z|}(xz)(yt)+(-1)^{|t|(|y|+|z|)}(xt)(yz).&\label{SJ2}
\eee

In particular, for any odd elements $x,y\in J_{\1}$ we have $xy=-yx$.  We will denote below the product of odd elements $x,y$ in a Jordan superalgebra $J$ as $[x,y]$.
\smallskip

{\bf Example 1.}
Let $A=A_{\0}+A_{\1}$ be an associative superalgebra. The new operation
\bes
a\circ b=\tfrac12(ab+(-1)^{|a||b|}ba)
\ees
defines a structure of a Jordan superalgebra on $A$. We will denote this Jordan superalgebra as $A^{(+)}$. 

\smallskip
{\bf Example 2.} A linear operator $*:A\rightarrow A$ on a superalgebra $A$ is called a {\em superinvolution} if it satisfies identities $(a^*)^*=a,\ (ab)^*=(-1)^{|a||b|}b^*a^*$ for arbitrary elements $a,b\in A_{\0}\cup A_{\1}$.  If $A$ is associative  then the subspace of symmetric elements $H(A,*)=\{a\in A\,|\, a^*=a\}$ is a Jordan subsuperalgebra of $A^{(+)}$.

\smallskip
{\bf Example 3.} Let $V=V_{\0}+V_{\1}$ be a $\Z/2\Z$-graded vector space over $F$ with a bilinear form $(v\,|\,w)$ on $V$ such that $(v\,|\,w)$ is symmetric on $V_{\0}$, skew-symmetric on $V_{\1},$ and $(V_{\0}\,|\,V_{\1})=(0)$.  The direct sum of vector spaces $J(V)=F\cdot 1+V= (F\cdot 1+V_{\0})+V_{\1}$ is a Jordan  superalgebra with respect to multiplication $v\cdot w=(v|w)1.$ We refer to this superalgebra as  {\em superalgebra  of a superform}  and denote $J(V).$

\smallskip
{\bf Example 4.} The $3$-dimensional Kaplansky superalgebra $K_3=Fe+(Fx+Fy)$  with the multiplication $e^2=e, \,ex=\tfrac12 x,\, ey=\tfrac12 y, [x,y]=e$ is simple and not unital.

\smallskip
{\bf Example 5.} The $1$-parametric family of 4-dimensional superalgebras $D_t(F)=(Fe_1+Fe_2)+(Fx+Fy)$, where $e_1,e_2$ are orthogonal even idempotents, $e_ix=\tfrac12 x,\,e_iy=\tfrac12 y,\,[x,y]=e_1+te_2,\, t\in F$. The superalgebra $D_t(F)$ is simple if $t\neq 0$.

\smallskip
{\bf Example 6.} V.\,Kac introduced the 10-dimensional simple superalgebra $K_{10}$ that is related (via the Tits-Kantor-Koecher construction) to the exceptional 40-dimensional Lie superalgebra.

\smallskip
{\bf Example 7.} Let $A$ be an associative commutative superalgebra equipped with a bilinear superskew-symmetric map $[,]: A\times A\rightarrow A, \, [A_{\bar i},A_{\bar j}]\subseteq A_{\overline{i+j}}$. A {\em Kantor double} is a direct sum of vector spaces
\bes
Kan(A,[,])=A+Av, \, |v|=1,
\ees
with the product
\bes
a(vb)=(-1)^{|a|}vab,\, (vb)a=vba, \, (va)(vb)=(-1)^{|a|}[a,b].
\ees
We say that $[,]$ is a {\em Jordan bracket} if $Kan(A,[,])$ is a Jordan superalgebra with respect to the grading $Kan(A,[,])_{\0}=A_{\0}+A_{\1}v,\, Kan(A,[,])_{\1}=A_{\1}+A_{\0}v.$
I.\,Kantor \cite{Kan1} showed that every Poisson bracket is a Jordan bracket.  Another example of a Jordan bracket is a {\em bracket of vector type}: let $(A,d)$ be an associative commutative superalgebra with an even derivation $d:A\rightarrow A$, then the bracket $[a,b]=a^db-ab^d$ is Jordan.  D.\,King and K.\,McCrimmon \cite{KingMC} found graded identities that determined the class of Jordan brackets. If the bracket $[,]$ is fixed, then we talk just about superalgebra $Kan(A)$.
\smallskip

N.\,Cantarini and V.\,Kac \cite{CanKac} showed that Jordan brackets are in 1-1 correspondence with {\em contact brackets} (see \cite{Kac_vdL}).

\smallskip
{\bf Example 8.} Let $A$ be an associative commutative superalgebra equipped with a Jordan bracket $[,]:A\otimes A\rightarrow A$ as in the Example 7.  Suppose further that $A$ has a 
$\Z/2\Z$-grading  $A=A_{(0)}+A_{(1)}$ that is compatible with the grading $A=A_{\0}+A_{\1}$ and, moreover, $[A_{(i)},A_{(j)}]\subseteq A_{(i+j)}$. Then the {\em Kantor doble} $Kan(A,[,])$ has the  
$\Z/2\Z$-grading
\bes
Kan(A,[,])&=&Kan(A,[,])_{(0)}+Kan(A,[,])_{(1)},\\
Kan(A,[,])&=&A_{(0)}+A_{(1)}v,\ Kan(A,[,])_{(1)}=A_{(1)}+A_{(0)}v.
\ees
We will refer to $A_{(0)}+A_{(1)}v$ as a {\em twisted Kantor double}.  A twisted Kantor double may be a simple Jordan superalgebra that is not isomorphic to any Kantor double.

For example,  let $A=\R[\sin t,\cos t],$   (see \cite{ZhSh}) where $\R$ is the field of real numbers,  $A_{(0)}=\R[\sin 2t,\cos 2t], \, A_{(1)}=(\sin t)A_{(0)}+(\cos t)A_{(0)},\ A=A_{(0)}+A_{(1)}$ is a $\Z/2\Z$-grading,  $[f(t),g(t)]=f'(t)g(t)-f(t)g'(t))$.
The odd part $A_{(1)}v$ of the twisted Kantor double is a projective module over the even part $A_{(0)}$, it is 2-generated, but not 1-generated. Similar examples over an arbitrary field of zero characteristic were constructed in  \cite{Zhel2}. 
 Moreover,  the examples where the odd part  is an $n$-generated module over the even part but can not be generated by less then $n$ elements were constructed in \cite {ZhZ}  for  abitrary $n\geq 2$ (over fields $\R$ and $\C$).  

\smallskip

{\bf Example 9.} 
Let $(A,d)$ be an associative commutative superalgebra with an even derivation $d:A\rightarrow A$.  Following  \cite{ChengKac},  C.\,Martinez and E.\,Zelmanov introduced a family  of Cheng-Kac Jordan superalgebras  $JCK(A,d)$ that are free $A$-modules of rank 8.  The superalgebras $JCK(F[t,t^{-1}],d/dt)$ are related  via the Tits-Kantor-Koecher construction to the exceptional superconformal algebras $CK_6$ (see \cite{GLS,MarZel}). 

Observe that when $A$ is equipped with a $\Z/2\Z$-grading competible with the derivation $d$,  the superalgebra $JCK(A,d)$ also  contain the twisted subsuperalgebra whose odd part is not a free $A$-module \cite{Zhel4}.   We will call it the {\em twisted Cheng-Kac superalgebra}.

\smallskip

Let $M_k(F)$ denote the algebra of $k\times k$ matrices over $F$., and let $k=m+n$. The grading
\bes
A=M_{m+n}(F)=A_{\0}+A_{\1},\ A_{\0}=\left(\begin{array}{cc}
*&0\\
0&*\end{array}\right),\ 
A_{\1}=\left(\begin{array}{cc}
0&*\\
*&0
\end{array}\right)
\ees
gives rise to a simple associative superalgebra.  

C.T.C.\,Wall \cite{Wall} proved that every finite dimensional associative superalgebra over an algebraically closed field $F$ is isomorphic to $M_{n+m}(F)$ or $M_n(\sqrt 1)$.

\smallskip
Associative superalgebras $M_{n+m}(F), \ M_n(\sqrt 1)$ give rise to simple Jordan superalgebras $M_{n+m}(F)^{(+)}$ and  $M_n(\sqrt 1)^{(+)}$ for $n>1$.

\smallskip
If $n=2k$ is even then the superalgebra $M_{m+2k}$ is equipped with the orthosymplectic involution $*$. The Jordan superalgebra of symmetric elements $H(M_{m+2k}(F),*)$ is called the {\em Jordan orthosymplectic superalgebra} and denoted $Josp_{n,2k}(F)$.
\smallskip

The associative superalgebra $M_{n+n}(F)$ is equipped with another superinvolution 
\bes
\left(\begin{array}{cc}
a&b\\
d&d\end{array}\right)\rightarrow \left(\begin{array}{cc}
d^t&-b^t\\
c^t&a^t\end{array}\right),
\ees
where $t$  is the transposition. The Jordan superalgebra of symmetric elements is denoted as $JP_n(F)$.
\smallskip

V.\,Kac \cite{Kac} (see also I.\,Kantor \cite{Kan1}) proved that every simple finite dimensional Jordan superalgebra over an algebraically closed field $F$ of zero characteristic is isomorphic to one of the superalgebras $M_{m+n}(F)^{(+)},$ $M_n(\sqrt 1)^{(+)},$   $Josp_{n,2k}(F)$, $JP_n(F)$,  a superalgebra of a superform,  $K_3,\, D_t,\, K_{10}$ or $Kan(G_n)$.

In \cite{RZ} it was shown that if $char\,F=p>3$ and the even part $J_{\0}$ is semisimple then every simple finite dimensional Jordan superalgebra is also isomorphic to one of the examples above.

If $char\,F=3$ then some new examples appear (see \cite{RZ, Sh3}).

If the even part $J_{\0}$ is not semisimple then the only  new examples are Kantor doubles 
$Kan(O_n\otimes G_m)$ \cite{MarZel,GLS} and Cheng-Kac superalgebras $JCK(O_n\otimes G_m,d)$, where $O_n=F[t_1,\ldots,t_n\,|\,t_i^p=0,\, 1\leq i\leq n]$ is the algebra of trancated polynomials, $d$ is an even derivation.

\smallskip

Now let us discuss infinite dimensional simple Jordan algebras and superalgebras.
In \cite {Zel4} it was proved that every simple Jordan algebra (finite or infinite dimensional) is isomorphic to $A^{(+)}$, where $A$ is a simple associative algebra,  or  to $H(A,*)$, where 
$A$ is a simple associative algebra with an involution $*:A\rightarrow A$, or to Jordan algebra of a nondegenerate symmetric bilinear form in a vector space over some extension of the ground field $F$, or to an exceptional Albert algebra that is 27-dimensional over its center.

Now we will formulate a slightly modified conjecture on classification of simple (finite or infinite dimensional) Jordan superalgebras that is due to N.\,Cantarini and V.\,Kac \cite{CanKac}.

\begin{conj} A simple Jordan superalgebra over a field $F$ (with a nonzero odd part) is isomorphic to one of the following superalgebras:
\begin{itemize}
\item [I)]
$A^{(+)}$, where $A$ is a simple associative superalgebra with a noncommutative even part;
\item[II)]
$H(A,*)$, where $A$ is a simple associative superalgebra with a superinvolution $*$ and with a noncommutative even part (see \cite{GA, GA-M});
\item[III)]
a Jordan superalgebra of a superform or $K_3$,  or $D_t$,  or $K_{10}$ over some extension of the ground field $F$;
\item[IV)]
a Kantor double  $Kan\,(A,[,])$ or a twisted Kantor double.
\item[V)]
a Jordan Cheng-Kac superalgebra $JCK(A,d)$, where $A$ is an associative commutative superalgebra, $d$ is an even derivation of $A$ and $A$ does not have  nontrivial $d$-invariant ideals (or a twisted subsuperalgebra of $JCK(A,d)$).
\end{itemize}
\end{conj}

N.\,Cantarini and V.\,Kac \cite{CanKac} proved the Conjecture for lineary compact simple Jordan superalgebras. V.\,Kac, C.\,Mart\'inez and E.\,Zelmanov \cite{KMZ} proved the Conjecture for superconformal  Jordan algebras, i.e. graded Jordan superalgebras $J=\sum_{i\in\Z}J_i$ that are  simple and have all dimensions $\dim_FJ_i$ uniformarly bounded.

\smallskip

A Jordan superalgebra $J$ is called {\em special} if it is embeddable in $A^{(+)}$ for some associative superalgebra $A$. A superalgebra is {\em $i$-special} if it is a homomorphic image of a special superalgebra. Equivalently, a superalgebra is $i$-special if it satisfies all identities that are satisfied by all special Jordan superalgebras.

It is easy to see that examples of types I -- III, exept for $K_{10}$, are special. Yu.\,Medvedev and E.\,Zelmanov \cite{MedZel1} proved that the superalgebra $K_{10}$ is not even $i$-special. 
K.\,McCrimmon \cite{McC}  (see also \cite{Sh5})  proved that a Kantor double $K(A,[,])$ of a bracket of vector type is special.  He also proved that the Kantor double of the algebra of the classical Poisson bracket is not special.  On the other hand,  I.\,Shestakov \cite{Sh4, MarShZ} showed that Kantor doubles are always $i$-special.   C.\,Martinez, I.\,Shestakov and E.\,Zelmanov \cite{MarShZ} showed that Jordan Cheng-Kac superalgebras are special.
\medskip

The purpose of this paper is a  partial progress towards classification of simple infinite dimensional Jordan superalgebras.  
The first theorem shows that we may restrict ourself to superalgebras with infinite dimensional even part.  
We consider the superalgebras whose even parts are infinite dimensional algebras of ``Clifford type'',  that is,   direct sums of algebras of bilinear form.   The results of \cite{RZ} show that the number of summonds in these sums is 1 or 2.   In the theorems 2 and 3 we consider the two cases.  It occurs that the only simple infinite dimensional Jordan superalgebras of this type are the superalgebras of superforms.
\\[2mm]
{\bf Theorem 1.}
{\em Let $J=J_{\0}+J_{\1}$ be a simple unital Jordan superalgebra with finite dimensional even part $J_{\0}.$ Then  either $J$ is finite dimensional or $J$ is the Jordan superalgebra of a nondegenerate supersymmetric bilinear form on a vector superspace $V=V_{\0}+V_{\1},\,  \dim_FV_{\0}<\infty,\, \dim_FV_{\1}=\infty$.
}\\[2mm]
{\bf Theorem 2.}
{\em Let $J=J_{\0}+J_{\1}$ be a simple unital Jordan superalgebra with $J_{\0}=J(V)$ beeing a Jordan algebra of a nondegenerate symmetric bilinear form. Then either $\dim_FJ<\infty$ or $J$ is a superalgebra of a nondegenerate supersymmetric bilinear superform.}\\[2mm]
{\bf Theorem 3.}
{\em Let $J=J_{\0}+J_{\1}$ be a simple unital Jordan superalgebra with $J_{\0}$ is isomorphic to the direct sum of two simple Jordan algebras of  nondegenerate symmetric bilinear forms. Then  $\dim_FJ<\infty$.}\\[2mm]
%{\bf Theorem 3.} {\em Let $J=J_{\0}+J_{\1}$ be a simple unital Jordan superalgebra over a field $F$ of zero characteristic. Suppose that  $J_{\0}$ contains a nonzero nilpotent ideal $N$ such that $J_{\0}/N$ is a Jordan algebra of a  nondegenerate symmetric bilinear form. Then  $J\cong Kan(G_n),\, n\geq 1$, the Kantor double of a Grassmann superalgebra.}
\smallskip

Throughout the paper,  without loss of generality we assume that the ground field $F$ is algebraically closed.

For a superalgebra $J=J_{\0}+J_{\1}$, we denote $A=J_{\0},\, M=J_{\1}$.  We will denote multiplication of even elements or an even element and an odd element  in $J$   by $\circ, \,\cdot$ or juxtoposition. The product of the odd elements $x,y\in M$  is denoted as $[x,y]$. 

On a Jordan (super)algebra $J$ we consider Jordan triple product
\bes
\{x,y,z\}=(xy)z+x(yz)-(-1)^{|y||z|}(xz)y,
\ees
for $x,y,z\in A\cup M. $ If $a\in A$ then denote 
\bes
U(a):J\rightarrow J,\ U(a):x\mapsto \{a,x,a\}.
\ees
For an element $a\in A\cup M$ let $R(a)$ denote the multiplication operator 
\bes
R(a):J\rightarrow J,\ x\mapsto xa.
\ees
The Jordan super-identity implies the following operator super-identity
\bee\label{Jordan_operator}
&R(x)R(y)R(z)+(-1)^{|x||y|+|x||z|+|y||z|}R(z)R(y)R(x)+(-1)^{|y||z|}R((xz)y)&\nonumber \\
&=R(xy)R(z)+(-1)^{|y||z|}R(xz)R(y)+(-1)^{|x||z|+|y||z|}R(yz)R(x).&
\eee
The operator $D(x,y)=R(x)R(y)-(-1)^{|x||y|}R(y)R(x)$ is an even or odd derivation depending on parity $|x|+|y|$.  We will often use the following identity for $a\in A,\,x\in A\cup M$:
\bes 
D(a^2,x)=2D(a,ax).
\ees
\section{Proof of Theorem 1}
\hspace{\parindent}
Let $S=S_{\0}+S_{\1}$ be a simple (not, necessary Jordan) superalgebra  over a field $F$.  We start with discussion of the centroid $Cent\,S$. This part may be folklorically known. We do it for the benefit of a reader.
\smallskip

Let $M(S)$ be the {\em multiplication algebra}  of $S$, that is, the  algebra of linear transformations of the vector space $M$ generated by all left and right multiplications by elements from $S_{\0}\cup S_{\1}$.  The algebra $M(S)$ is also $\Z_2$-graded: 
$$
M(S)=M(S)_{\0}+M(S)_{\1}. 
$$

The centroid $Cent\,(S)$ is the subalgebra of all even linear transformations of $S$ that commute with all elements form $M(S)$. 
Standard arguments show that the centroid of a simple superalgebra is a field extension of $F$.  We call the superalgebra $S$ {\em central} if $Cent\,(S)=F$.

For $i=0$ or $1$ let $Cent\,(S_{\bar i})$ be the algebra of linear transformations of the space $S_{\bar i}$ that commute with all elements from the restriction of $M(S)_{\0}$ to $S_{\bar i}$. For an arbitrary transformation $\f$ from $Cent\,(S)$  the restriction of $\f$  to $S_{\bar i}$ lies in $Cent\,(S_{\bar i})$. We denote this restriction as $\pi_i(\f)$.
\begin{lem}\label{lem1.1}
$\pi_i: Cent\,(S)\rightarrow Cent(S_{\bar i}),\ i=0,1$  is an isomorphism.
\end{lem}
\prf
We need to show that the mapping $\pi_i$ is surjective.  Let $\psi\in Cent\,(S_{\bar i})$ and let $a\in S_{\overline{1-i}}$. Since the superalgebra $S$ is simple,  there exist elements $a_k\in S_{\bar i}$ and operators $W_k\in M(S)_{\1}$ such  that 
$$ 
a=\sum_k a_kW_k.
$$
We will define the extension $\tilde\psi$ of $\psi$ to $S_{\overline{1-i}}$ via
$$
\tilde\psi(a)=\sum_k\psi(a_k)W_k
$$
It remains to prove that $\sum_k a_kW_k=0$ implies $\sum_k \psi(a_k)W_k=0$.  If $\sum_k \psi(a_k)W_k\neq 0$ then 
$$
(\sum_k \psi(a_k)W_k)M(S)_{\1}=S_{\0},
$$
because of simplicity of the superalgebra $S$.

In particular, there exists an operator $U\in M(S)_{\1} $ such that 
$$
(\sum_k \psi(a_k)W_k)U\neq 0.
$$
Since $W_kU\in M(S)_{\0}$ and $\psi\in Cent\,(S_{\bar i})$,  it follows that 
$$
\sum_k \psi(a_k)(W_kU)=\psi(\sum_k a_kW_kU)= 0,
$$
a contradiction. 

The construction above implies that for an arbitrary element $b\in S_{\bar i}$  and an arbitrary operator $W\in M(S)_{\1}$ we have 
$$
\tilde\psi(bW)=\psi(b)W.
$$
Now let $b\in S_{\overline{1-i}},\, W\in M(S)_{\0}$. We need to show that $\tilde\psi(bW)=\tilde\psi(b)W$.
Let $b=\sum_ka_kW_k,\, a_k\in S_{\bar i},\, W_k\in M(S)_{\1}$. Then
\bes
\tilde\psi(bW)=\tilde\psi(\sum_k a_kW_kW)=\sum_k\psi(a_k)W_kW=(\sum_k\tilde\psi(a_kW_k))W=\tilde\psi(b)W.
\ees
Finally, let $b\in S_{\overline{1-i}},\, W\in M(S)_{\1},\, b=\sum_ka_kW_k,\,a_k\in S_{\bar i},\, W_k\in M(S)_{\1}.$ Then 
$bW=\sum_k a_kW_kW,\, W_kW\in M(S)_{\0}$. Hence 
$$
\tilde\psi(b)=\sum_k\psi(a_k)W_k,\,\ \psi(bW)=\sum_k \psi(a_k)W_kW=\tilde\psi(b)W.
$$
We proved that the mapping $\pi_i$ is injective,  hence it is an isomorphism. This completes the proof of the lemma.

\ctd

\begin{lem}\label{lem1.2}
Let $S$ be a central simple superalgebra over a field $F$. Suppose that the multiplication algebra $M(S)$ satisfies a polynomial identity. Then the superalgebra $S$ is finite dimensional.
\end{lem}
\prf
Each homogeneous component $S_{\bar i}, \, i=0,1,$ is an irreducible module over $M(S)_{\0}$.  By Lemma \ref{lem1.1} and centrality of the superalgebra $S$,  the centralizer of this module is the ground field $F$. The celebrated  Kaplansky Theorem \cite[chapter 6]{Her} implies that if $A$ is an associative $PI$-algebra and $M$ is an irreducible right $A$-module with  the centralizing division algebra $\Delta$ then $M$ is finite dimensional over $\Delta$.  Hence,  the vector space $S_{\bar i}$ is finite dimensional.  This completes the proof of the lemma.

\ctd

\smallskip

From now on in this section $J=A+M$ is a simple unital Jordan superalgebra over a field $F$, $\dim_FA\leq n$.  Without loss of generality we will assume that $F=Cent\,(J)$. Moreover, if $\tilde F$ is an algebraically closed field extension of $F$ and $card\,\tilde F>\dim_FJ$ then the superalgebra $\tilde J=\tilde F\otimes_FJ$ is simple, $\dim_{\tilde F}\tilde J_{\0}=\dim_FA$.  Therefore without loss of generality  we assume that the field $F$ is algebraically closed and $card\,F>\dim_F J$.

\begin{lem}\label{lem1.3}
For an arbitrary finitely generated subsuperalgebra $J'\subseteq J$ there exist $m\geq 1$ and homomorphisms $\f_i:J'\rightarrow \bar J_i,\,1\leq i\leq m,$ into finite dimensional simple Jordan superalgebras $\bar J_i$,  such that $\f_1\oplus\cdots\oplus\f_m:J'\rightarrow \bar J_1\oplus\cdots\oplus \bar J_m$ is an embedding,  and dimensions of all even parts $(\bar J_i)_{\0},\, 1\leq i\leq m, $ are $\leq n$.
\end{lem}
\prf
Every finitely generated $A$-subbimodule of $M$ is finite dimensional (see \cite{Jac}).  If $M'$ is an $A$-subbimodule of $M$ then $A+M'$ is a subsuperalgebra of $J$.  This implies that the superalgebra $J$ is locally finite dimensional, i.e. every finitely generated subsuperalgebra of $J$ is finite dimensional.

We will show that for an arbitrary nonzero element $a\in A\cup M$ and an arbitrary finite dimensional subsuperalgebra $J',\ a\in J'\subseteq J,$ there exists a homomorphism $\f:J'\rightarrow \bar J$ into a simple finite dimensional Jordan superalgebra $\bar J$ with $\dim(\bar J)_{\0}\leq n$ such that $\f(a)\neq 0$.

Indeed, there exists a multiplication operator $W\in M(J)$ such that $aW=1$. Let $J''$ be a finite dimensional  subsuperalgebra of $J$  that contains $1, J'$ and all  elements involved in the operator $W$.  Since $J''\ni 1, $ it follows that exists an epimorphism $\psi:J''\rightarrow H$ onto a simple finite dimensional superalgebra $H,\, \psi(1)=1$.  Clearly,  $\dim_F H_{\0}\leq\dim_F A=n$. The restriction of $\psi$ to $J'$ maps $a$ to a nonzero element.

Now let $J'$ be an arbitrary finite dimensional subsuperalgebra of $J$. Let $\f_i:J'\rightarrow \bar J_i,\,1\leq i\leq m$,  be homomorphisms    into simple finite dimensional superalgebras, $\dim_F(\bar J_i)_{\0}\leq n$, such that $\cap_{i=1}^m\ker \f_i$ has minimal dimension.
In view of the above $\cap_{i=1}^m\ker \f_i=(0)$.  This completes the proof of the lemma.

\ctd

The following lemma immediately follows from classification of  simple finite dimensional Jordan superalgebras.
\begin{lem}\label{lem2.0}
There exists a  a function $f:\N\rightarrow \N$ with the following property: 
 if $J=J_{\0}+J_{\1}$ is a simple finite dimensional Jordan superalgebra and $\dim_FJ_{\0}\leq n$ then $\dim_FJ\leq f(n)$ or $J$ is a superalgebra of a superform.
\end{lem}
\prf
An arbitrary simple finite dimensional Jordan superalgebra over an algebraically closed field $F$ of characteristic $\neq 2$ either has dimension $\leq 10$ or belongs to one of the following families:
\bes
&M_{m+n}(F)^{(+)}, \, M_n(\sqrt 1)^{(+)},\, Josp_{n,2k}(F),\, JP_n(F),&\\
&Kan(O_n\otimes G_m),\, JCK(O_n\otimes G_m,d) \hbox{  (see \cite{Kac,  Kan1,  MarZel,  RZ} or the Introduction)}.&
\ees
For a superalgebra $J=J_{\0}+J_{\1}$ from one of the families $M_n(\sqrt 1)^{(+)},\, JP_n(F)$, $  Kan (O_n\otimes G_m), \, JCK(O_n\otimes G_m,d)$ we have $\dim_FJ_{\0}=\dim_FJ_{\1}$.
For a superalgebra $J=M_{m+n}(F)^{(+)}$ we have $\dim_FJ_{\1}\leq\dim_FJ_{\0}$.  Finally, for a superalgebra $J=Josp_{n,2k}(F)$ we have $\dim_FJ_{\1}\leq 2\dim_FJ_{\0}$.
In any case, $\dim_FJ_{\1}\leq f(\dim_FJ_{\0})$, where $f(n)=\max(2n,10)$.  This completes the proof of the lemma.

\ctd

{\bf\underline{Proof of Theorem 1.}}\\[2mm]

Let $J=A+M$ be a simple unital Jordan superalgebra, $\dim_FA=n$.

If the multiplication algebra $M(J)$ is $PI$ then the superalgebra $J$ is finite dimensional by Lemma \ref{lem1.2}. Suppose that the algebra $M(J)$ is not $PI$. 

\smallskip
Let $d=2f(n)$. There exist $d+1$ operators $W(a_1,a_2\ldots),\,W_i(a_{i1},a_{i2},\ldots),\,1\leq i\leq d,\, a_i, a_{ij}\in J_{\0}\cup J_{\1}$ and an element $a\in A\cup M$ such that 
$$
1=aS_d(W_1,\ldots,W_d)W,
$$  where 
$$
S_d=\sum_{\s\in\Sigma_n} (-1)^{sgn\,\s} x_{\s(1)}\cdots x_{\s(n)}
$$
is the standard noncommutative polynomial.  

Let $J'$ be an arbitrary  finite dimensional subsuperalgebra of $J$.  We want to prove that  $J'$  is embeddable into a direct sum of Jordan superalgebras of superforms.
Without loss of generality, we may assume that $J'$  contains $1,a,$ and all the elements $a_{ij}, \, a_i$.  Clearly, $\dim J'_{\0}\leq n$. Let $\f:J'\rightarrow\bar J$
be a homomorphism into a simple Jordan superalgebra $\bar J$ of dimension $\leq f(n)$.  By the Amitsur-Levitzki Theorem, the multiplication algebra $M(\bar J)$ satisfies the identity $S_d(x_1,\ldots,x_d)=0$,  hence $\f(1)=\f(aS_d(W_1,\ldots,W_d)W)=0$, a contradiction.  Therefore,  $\dim\bar J>f(n)$ and by lemma \ref{lem2.0} $\bar J$ is a superalgebra of superform.

In a Jordan superalgebra of a  superform for arbitrary odd elements $x,y$ the operator $R([x,y])$ lies in the centroid. 
 Hence the same is true for the superalgebra $J$.  Since the superalgebra $J$ is central we conclude that 
 $$
 [M,M]\subseteq F\cdot 1.
 $$ 
 Let $Rad\,A$ be the radical of the finite dimensional Jordan algebra $A$. Let $J'$ be a finite dimensional subsuperalgebra of $J, \,J'_{\0}=A$, and let $\f:J'\rightarrow \bar J$ be one of the homomorphisms that we discussed above, $\bar J$ is a simple Jordan superalgebra of a superform.  The homomorphism $\f$ is surjective on even part, $\f(A)=\bar J_{\0} $,  hence $\f(Rad\,A)\subseteq Rad\,(\bar J_{\0})=(0)$.
 This implies that the algebra $A$ is semisimple.  Let $A=A_1\oplus\cdots\oplus A_m$  be the decomposition of $A$ into a direct sum of simple algebras, and let $e_i$ be the identity of the algebra   $A_i$.  
 
 Suppose that $m\geq 2$. Then 
 $$
 [\{e_i,M,e_i\},\{e_i,M,e_i\}]\subseteq Fe_i\cap F\cdot 1=(0),
 $$
 and for arbitrary $j,\,j\neq i,$
 $$
 [\{e_i,M,e_i\},\{e_i,M,e_j\}]\subseteq \{e_i,A,e_j\}=(0).
 $$
 This implies that $\{e_i,M,e_i\}$ is an ideal in $J$, hence 
 $$
 \{e_i,M,e_i\}=0,\ M=\oplus_{i\neq j}\{e_i,M,e_j\}.
 $$
 If $m\geq 3$ then arguing as above we conclude that $\{e_i,M,e_j\}$ is an ideal in $J$,  hence $\{e_i,M,e_j\}=(0)$ and $M=(0)$.
 In particular, $J=A$ is finite dimensional.
 
 Hence $A=A_1\oplus A_2,\, M=\{e_1,M,e_2\}$.  Let $(0)\neq [M,M]\subseteq F\cdot 1$.
 
 Let $J'$ be a finite dimensional subsuperalgebra of $J,\ J'_{\0}=A,\,  [J'_{\1},J'_{\1}]\neq (0)$,  and 
 let $\f:J'\rightarrow \bar J$ be a homomorphism into a simple Jordan superalgebra of a superform,  surjective on the even part. 
 If one of the ideals $A_1,A_2$ lies in $\ker\f$ then $J'_{\1}=\{e_1,J'_{\1},e_2\}\subseteq\ker\f$,  hence $1\in \ker\f$,  a contradiction.
 
 If $\f:A_1\oplus A_2\rightarrow\bar J_{\0}$ is an isomorphism then $A_i=Fe_i,\, i=1,2$.  In this case $J$ is isomorphic to a superalgebra of  a superform on a space $V=V_{\0}+V_{\1},\, \dim V_{\0}=1$.
 
 \smallskip
 Now assume that $m=1$, that is, the algebra $A$ is simple.  Again choose a finite dimensional subsuperalgebra $J'\subseteq J$ such that $J'_{\0}=A,\, [J'_{\1},J'_{\1}]\neq (0)$, and a homomorphism $\f:J'\rightarrow \bar J$, where $\bar J$ is a simple Jordan superalgebra of a superform,  $\f:A\rightarrow \bar J_{\0}$ is an isomorphism.  Let $A_0$ be the subspace of elements of zero trace, $A_0=Span\la (ab)c-a(bc)\,|\,a,b,c\in A\ra$.  Then $\f(A_0)\cdot \bar J_{\1}=(0)$,  which implies that $A_0\cdot M=(0)$.  Hence $J$ is a superalgebra of a superform on the vector space $A_0+M$.   This completes the proof of Theorem 1.
 
 \ctd

\section{Proof of Theorem 2}
\hspace{\parindent}
Let $J=A+M$ be a simple Jordan superalgebra, $A=F\cdot 1+V,$ where $V$ is an infinite dimensional vector space over $F$. 

Choose $2n$ orthonormal elements $v_1,\ldots,v_{2n}\in V, v_i\circ v_j=\delta_{ij},\, 1\leq i,j\leq 2n$.  Let $Cl(V)$ be the Clifford algebra of the nondegenerate bilinear form on $V$. Both the superalgebra $J$ and the Clifford algebra $Cl(V)$ are graded by the finite abelian group $(\Z/2\Z)^{2n}$. Let us define the gradings.

\smallskip
The operators $U(v_i),\,1\leq i\leq 2n,$ are automorphisms of the superalgebra $J$. They commute and $(U(v_i))^2=1$. Hence these operators generate the group $(\Z/2\Z)^{2n}\leq Aut\,J$. The superalgebra $J$ decomposes as a sum of eigenspaces
$$
J=\oplus_{\a\in (\Z/2\Z)^{2n}}J_{\a}.
$$
An element  $a$ lies in $J_{\a}$ for $\a=(\a_1,\ldots,\a_{2n}),\, \a_k=0$ or 1,  if and only if  $U(v_i)(a)=(-1)^{\a_i}a$ for any $i, 1\leq i\leq 2n$.

The operators $Cl(V)\rightarrow Cl(V),\, x\mapsto v_ixv_i,\,x\in Cl(V),$ are automorphisms of the algebra $Cl(V),\, 1\leq i\leq 2n$. They commute and generate the subgroup  $(\Z/2\Z)^{2n}\leq Aut\,Cl(V)$. The algebra $Cl(V)$ decomposes as a sum of eigenspaces
$$
Cl(V)=\oplus_{\a\in (\Z/2\Z)^{2n}}Cl(V)_{\a}.
$$
Again, $a$ lies in $Cl(V)_{\a}$ for $\a=(\a_1,\ldots,\a_{2n})$ if and only if  $v_iav_i=(-1)^{\a_i}a,\,1\leq i\leq 2n$.
\begin{lem}\label{lem2.1}
For an arbitrary $2n$-tuple $\a\in (\Z/2\Z)^{2n}$ there exists a unique element $v_{i_1}\cdots v_{i_k},$ $1\leq i_1<\cdots <i_k\leq 2n,$  where the product is considered in the Clifford algebra $Cl(V)$, such that $v_{i_1}\cdots v_{i_k}\in Cl(V)_{\a}$.
\end{lem}
\prf
Let $\a=(\a_1,\ldots,\a_{2n}),\, \a_i=0$ or 1, $1\leq i\leq 2n.$
Let $X=\{i\,|\,\a_i=0\},\, Y=\{i\,|\,\a_i=1\},\,|X|+|Y|=2n.$
If the orders $|X|,|Y|$ are even and $Y=\{j_1<\cdots <j_{2r}\}$ then the element $v_{j_1}\cdots v_{j_{2r}}$ lies in $Cl(V)_{\a}$.  If the orders $|X|,|Y|$ are odd and $X=\{i_1<\cdots <i_{2t+1}\}$ then the element $v_{i_1}\cdots v_{i_{2t+1}}$ lies in $Cl(V)_{\a}$.  

Let us prove uniqueness. Suppose that products $v_{i_1}\cdots v_{i_p},\, v_{j_1}\cdots v_{j_q}$
lie in the  same eigenspace $Cl(V)_{\a},\, i_1<\cdots<i_p,\, j_1<\cdots<j_q$ and $\{i_1,\ldots,i_p\}\neq \{j_1,\ldots,j_q\}$. Then the element  $v_{i_1}\cdots v_{i_p}v_{j_1}\cdots v_{j_q}$ belongs to the eigenfunction $(0,\ldots,0)\in (\Z/2\Z)^{2n}$. In other words, 
$$
v_i(v_{i_1}\cdots v_{i_p}v_{j_1}\cdots v_{j_q})v_i=v_{i_1}\cdots v_{i_p}v_{j_1}\cdots v_{j_q}
$$
for any $1\leq i\leq 2n$. 
Since $\{i_1,\ldots,i_p\}\neq\{j_1,\ldots,j_q\}$, it follows that 
$$
v_{i_1}\cdots v_{i_p}v_{j_1}\cdots v_{j_q}=\pm v_{k_1},\cdots v_{k_r},\, 1\leq k_1<\cdots<k_r\leq 2n,\, r\geq 1.
$$
If $r$ is even then
$$
U(v_{k_1})(v_{k_1}\cdots v_{k_r})=-v_{k_1}\cdots v_{k_r}
$$
If $r$ is odd then for any $j\in\{1,2,\ldots,2n\}\setminus\{k_1,\ldots,k_r\}$ we have
$$
U(v_j)(v_{k_1}\cdots v_{k_r})=-v_{k_1}\cdots v_{k_r}.
$$
Both cases contradict our assumption. This completes the proof of the lemma.

\ctd

If elements $a\in J_{\a}$ and $v_{i_1}\cdots v_{i_k}\in Cl(V)_{\a}$, $1\leq i_1<\cdots <i_k\leq 2n,$ belong to the same eigenvalue $\a$, then we write $a\sim v_{i_1}\cdots v_{i_k}$ or $a=a(i_1\cdots i_k)$.

If $a\sim v_{i_1}\cdots v_{i_k},$ $b\sim v_{j_1}\cdots v_{j_q},$  then both $a\circ b,\, a\cdot b\sim  v_{i_1}\cdots v_{i_k}v_{j_1}\cdots v_{j_q}$,  where $a\circ b$ means the product in $J$ and $a\cdot b$ stands for the product in $Cl(V)$.

The product of $a$ and $b$ may be equal to zero, but zero belongs to all eigenvalues.

Let $0\neq v\in V$ be an eigenvector. Then $v\sim v_i,\,1\leq i\leq 2n$ if $v\in Fv_i$ or $v\sim v_1\cdots v_{2n}$ if $v\in\{v_1,\ldots,v_{2n}\}^{\perp}=\{u\in V\,|\,u\circ v_i=0,\,1\leq i\leq 2n\}$.
\begin{lem}\label{lem2.2}
Let $ x=x(1\,2\ldots 2p)\in M,\,1\leq p<n.$ Then $[x,M]\subseteq \{v_1,\ldots,v_{2n}\}^{\perp}$.
\end{lem}
\prf
If $y\in M$ is an eigenvector and $[x,y]\neq 0$ then there are the following four options for the eigenfunction of $y$:
\begin{enumerate}
\item 
$y=y(1\,2\ldots 2p)$. In this case we assume that $[x,y]=1$;
\item
$y=y(1\,2\ldots \hat i \ldots 2p)$, we assume $[x,y]=v_i$;
\item
$y=y(1\,2\ldots 2p\,j), 2p<j,$ we assume $[x,y]=v_j$;
\item
$y=y(2p+1\ldots\,2n),\, [x,y]=u\neq 0,\, u\in \{v_1,\ldots,v_{2n}\}^{\perp}$.
\end{enumerate}
{\bf\underline{Option 1.}} \  
Observe that  $v_{2p}\circ x=D(x,v_{2p+1})=0.$ 
In fact,  we have $U(v_{2p})(x)=-x$, hence $(x\circ v_{2p})\circ v_{2p}=0$.  Since $R(v_{2p})^3=R(v_{2p})$, we have $x\circ v_{2p}=0$.  Furthermore,  $U(v_{2p+1})(x)=x$,
hence $x=(x\circ v_{2p+1})\circ v_{2p+1}$, and we have 
$$
D(x,v_{2p+1})=D((x\circ v_{2p+1})\circ v_{2p+1},v_{2p+1})=\tfrac12 D(x\circ v_{2p+1},1)=0.
$$
Consider now  the eigenvector
$$
y'(1\,2\ldots 2p-1\,2p+1)=(y\circ v_{2p+1})\circ v_{2p}.
$$
By the previous consideration,
\bes
v_{2p}D(x,y')&=&[v_{2p}\circ y',x]=[y\circ v_{2p+1},x]=[y,x]\circ v_{2p+1}+yD(v_{2p+1},x)\\
&=&[y,x]\circ v_{2p+1}=-v_{2p+1}.
\ees
Consider the nonzero element (we assume that $x\neq 0$)
\bes
x'&=&x'(1\,2\ldots 2p\, 2p+1)=x\circ v_{2p+1}=-x\circ (v_{2p}D(x,y'))\\
&=&-(x\circ v_{2p})D(x,y')+(xD(x,y'))\circ v_{2p}.
\ees
We have $x\circ v_{2p}=0$,  $xD(x,y')=[x(1\,2\ldots 2p),y'(1\,2\ldots 2p-1\,2p+1)]\circ x$,
and $A\ni [x,y']\sim v_{2p}v_{2p+1}$.  Hence $[x,y']=0$ and $x'=0$, a contradiction.

Therefore, the option 1 is impossible.

\smallskip
{\bf\underline{Option 2.}} \  Let $y=y(1\,2\ldots 2p-1),\  [x,y]=v_{2p}$. \\[2mm]
Let $y'=y'(1\,2\ldots 2p-2)=y(1\,2\ldots 2p-1)\circ v_{2p-1}$.Then $y=y'\circ v_{2p-1}$.
We have 
\bes
v_{2p-1}D(x,y')=[v_{2p-1}\circ y',x]=[y,x]=-v_{2p}.
\ees
{\bf\underline {Case 2.1.}}\ $n\geq 2p$.\\
Then $2n-2p+1\geq 2p+1$. Consider the element
\bes 
x'(2p+1\ldots 4p)=x(1\,2\ldots 2p)D(v_1,v_{2p+1})\cdots D(v_{2p},v_{4p})\neq 0,
\ees
and the element $x''=x'\circ v_{2p},\, x''=x''(2p\,2p+1\cdots 4p)$. We have 
\bes
x''=-x'\circ(v_{2p-1}D(x,y'))=-(x' \circ v_{2p-1})D(x,y')+(x'D(x,y'))\circ v_{2p-1}.
\ees
Let $x'''=x'\circ v_{2p-1}=x'''(2p-1\,2p+1\ldots 4p)$. Now,
$$
[x'''(2p-1\,2p+1\ldots 4p),x(1\,2\ldots 2p)]\sim v_1\cdots \hat v_{2p-1}\cdots v_{4p}.
$$
 Hence $[x''',x]=0$.

Furthermore,
$$
[x'''(2p-1\,2p+1\ldots 4p),y'(1\,2\ldots 2p-2)]\sim v_1\cdots \hat v_{2p}\cdots v_{4p}.
$$
 Hence $[x''',y']=0$  and $(x'\circ v_{2p-1})D(x,y')=0$.
 
 Consider another summond $(x'D(x,y'))\circ v_{2p-1}$. We have
 \bes
 [x'(2p+1\ldots 4p),x(1\,2\ldots 2p)]&\sim &v_1v_2\cdots v_{4p},\\ \ 
 [x',x]\circ y'&\sim& v_{2p-1}v_{2p}\ldots v_{4p}.
 \ees
 Hence $([x',x]\circ y')\circ v_{2p-1}=0$.
\smallskip

Also
$$
[x'(2p+1\ldots 4p),y'(1\,2\ldots 2p-2)]\sim v_1\cdots\hat v_{2p-1}\hat v_{2p}\cdots v_{4p},
$$
hence $[x',y']=0$.

We proved that the Case 2.1. is impossible.\\[2mm]
{\bf\underline {Case 2.2.}}\ $n<2p$.\\
Choose a subset $S\subseteq\{1,2,\ldots,2p-1\}$ of order $4p-2n$ such that $S\ni 2p-1$. Since $2n-2p+1+|S|=2p+1$, there exists a nonzero element
$$
x'=x'(S,2p, 2p+1\ldots, 2n).
$$
As above, let $x''=x'\circ v_{2p}=x''(S,2p+1,\ldots,2n)\neq 0.$ We have
$$ 
	x'=x''\circ v_{2p}=-x''\circ (v_{2p-1}D(x,y')).
$$
Since $|S|$ is even and $S\ni 2p-1$, it follows that $x''\circ v_{2p-1}=0$.  Hence
\bes
x'&=&(x''D(x,y'))\circ v_{2p-1},\\
x''D(x,y')&=&[x'',x]\circ y'+[x'',y']\circ x.
\ees
Now, 
$$
[x'',x]=[x''(S,2p+1,\ldots,2n),x(1\,2\ldots 2p)].
$$
Let $v_S$ be the product of all $v_i,\,i\in S$. Then
$$
[x'',x]\sim v_1\cdots v_{2n}\cdot v_S.
$$
The length of the product on the right hand  side is $2n-|S|=4n-4p=4(n-p)$. This implies that $[x'',x]=0$.

Consider another summond:
\bes
[x'',y']&=&[x''(S,2p+1,\ldots,2n),y'(1\,2\ldots 2p-2)]\\
&\sim &v_1\cdots v_{2p-2}\cdot v_{2p+1}\cdots v_{2n}\cdot v_S.
\ees
The set $S$ contains $2p-1$ and does not contain $2p$. Hence the length of the product is equal to $2n-|S|=4(n-p)$ and again $[x'',y']=0$. 

We proved that the Case 2.2. is impossible,  and therefore the option 2 is impossible.\\[2mm]
{\bf\underline {Option 3.}}  Let $y'(1\,2\ldots 2p)=y(1\,2\ldots 2p\,j)\circ v_j$. Then $y=y'\circ v_j$. The inner derivation $D(x,v_j)=0$, hence 
\bes
[x,y]=[x,y'\circ v_j]=[x(1\,2\ldots  2p),y'(1\,2\ldots 2p)]\circ v_j=\hbox{ (by option 1) }=0.
\ees
This leaves us with the only option 4, which completes the proof of the lemma.

\ctd
\begin{lem}\label{lem2.3} 
Let $x=x(1\,2\ldots 2p+1)\in M,\,1\leq p<n$. Then $[x,M]\subseteq \{v_1,\ldots,v_{2n}\}^{\perp}.$
\end{lem}\prf
If $y\in M$ is an eigenvector and $[x,y]\neq 0$ then there are the following options for the eigenfunction of $y$:
\begin{itemize}
\item[5.]
$y=y(1\,2\ldots 2p+1)$;
\item[6.]
$y=y(1\,2\ldots \hat i\ldots 2p+1),\,1\leq i\leq 2p+1$;
\item[7.]
$y=y(1\,2\ldots 2p+1\,j),\, 2p+1<j\leq 2n$;
\item[8.]
$y=y(2p+2\ldots 2n)$.
\end{itemize}
{\bf\underline {Option 5.}} Let $y'(1\,2\ldots 2p)=y(1\,2\ldots 2p+1)\circ v_{2p+1}$. Then $y=y'\circ v_{2p+1}$. We have $D(x(1\,2\ldots  2p+1),v_{2p+1})=0$.  Hence
\bes
[x,y]=[x,y'\circ v_{2p+1}]=[x,y']\circ v_{2p+1}=\hbox{ (by lemma \ref{lem2.2})}= 0.
\ees
{\bf\underline {Option 6.}} $[x,y(1\,2\ldots \hat i\ldots 2p)]=0$ because on one hand $[x,y]\sim v_i$. On the other hand,  by lemma  \ref{lem2.2} $[x,y]\in\{v_1,\ldots,v_{2n}\}^{\perp}$.\\[2mm]
{\bf\underline {Option 7.}} Suppose that $p+1<n$.  Then on one hand $[x,y]\sim v_j$. On the other hand, by lemma  \ref{lem2.2} $[x,y]\in\{v_1,\ldots,v_{2n}\}^{\perp}$. Hence $[x,y]=0$.
\smallskip

Now suppose that $p=n-1,\,x=x(1\,2\ldots 2n-1),\,y=y(1\,2\ldots 2n),\, [x,y]=v_{2n}$.
Let $x'(1\,2\ldots 2n-2)=x\circ v_{2n-1},\, x=x'\circ v_{2n-1}$. Then $v_{2n}=[x,y]=[x'\circ v_{2n+1},y]=v_{2n-1}D(x',y)$ since $v_{2n-1}\circ y=0$.

There exists a nonzero element $x''(2\,3\ldots 2n-1\,2n)$.  Let $x'''=x''\circ v_{2n}$.  Then $x''=x'''\circ v_{2n}=x'''(v_{2n-1}D(x',y))=-x'''D(x',y)\circ v_{2n-1}$, since $x'''(2\,3\ldots 2n-1)\circ v_{2n-1}=0$.

We have 
$$
[x'''(2\,3\ldots 2n-1),x'(1\,2\ldots 2n-2)]\sim v_1v_{2n-1}, 
$$
hence $[x''',x']=0$.   Similarly,
$$
[x'''(2\,3\ldots 2n-1),y(1\,2\ldots 2n)]\sim v_1v_{2n}, 
$$
hence $[x''',y]=0$.  Therefore $x''=0$,  a contradiction.

Thus the option 7 is impossible. This completes the proof of the lemma.

\ctd

Denote the orthonormal system $v_1,\ldots,v_{2n}$ as $\V$. Let $M(\V,r),\,1\leq r\leq n,$ be the linear span of all elements in $M$ that are equivalent to $v_{i_1}\ldots v_{i_{2r}},\,1\leq i_1<\cdots<i_{2r}\leq 2n$,  or to $v_{i_1}\ldots v_{i_{2r+1}},\,1\leq i_1<\cdots<i_{2r+1}\leq 2n$.  Clearly, $M(\V,r)$ is a $(F\cdot 1+\sum_{i=1}^{2n} Fv_i)$-submodule of $M$.
\smallskip

Lemmas \ref{lem2.2},  \ref{lem2.3} imply that for $1\leq r<n$
\bes \ 
[M(\V,r),M]\subseteq \{v_1,\ldots,v_{2n}\}^{\perp}.
\ees
Let $\tilde\V=\{v_1,\ldots,v_{2m}\}$ be an orthonormal extension of $\V$, $n<m$.
\begin{lem}\label{lem2.4}
Let $1\leq r<n$. Then 
\bes
M(\V,r)\subseteq \sum_{t=1}^{m-1}M(\tilde \V,t).
\ees
\end{lem}
\prf
Let $x\in M,\,x\sim v_1\cdots v_{2r}$ with respect to $\V,\,1\leq r<n$.  Let
\bes
x=\sum_{\pi}\a_{\pi}x_{\pi}+\sum_{\tau}\b_{\tau}x_{\tau},
\ees
where $\a_{\pi},\b_{\tau}\in F;\,\pi,\tau$ are subsets of $\tilde\V,\, |\pi|=2t,\,|\tau|=2t+1,\,0\leq t\leq m$.
If $\pi=\{j_1<\cdots<j_{2t}\}$ then $x_{\pi}\sim v_{j_1}\cdots v_{j_{2t}}$,   if $\tau=\{j_1<\cdots<j_{2t+1}\}$ then $x_{\tau}\sim v_{j_1}\cdots v_{j_{2t+1}}$ with respect to $\tilde\V$.

The element $x$ belongs to the eigenvalue 1 with respect to $U(v_{2r+1})$,$\ldots,$  $U(v_{2n})$ and to the eigenvalue $-1$ with respect to $U(v_1),\ldots,U(v_{2r})$. Hence every subset $\pi$ contains $1,2,\ldots,2r$ and does not contain $2r+1,\ldots,2n$. Therefore,
\bes
2r\leq|\pi|\leq 2m-(2n-2r)\leq 2m-2.
\ees
Similarly,  every subset $\tau$ contains $2r+1,\ldots,2n$ and does not contain $1,2,\ldots,2r$.  Hence
\bes
2n-2r\leq|\tau|\leq 2m-2r.
\ees
We proved that $2\leq|\pi|, |\tau|\leq 2m-2$. Hence $x\in\sum_{t=1}^{m-1}M(\tilde\V,t).$

\smallskip
If $x\sim v_1\cdots v_{2r+1}$, then $x=x'\circ v_{2r+1}$,  where $x'\sim v_1\ldots v_{2r}$. 
We proved that $x'\in\sum_{t=1}^{m-1}M(\tilde\V,t)$.  Since $M(\tilde\V,t)\circ v_{2r+1}\subseteq \sum_{t=1}^{m-1}M(\tilde\V,t),$
this  completes the proof of the lemma.

\ctd 

\begin{lem}\label{lem2.5} Let $1\leq r<n$.Then $M(\V,r)=(0)$.
\end{lem}
\prf
Let $W$ be a finite-dimensional subspace of $V$. There exists an orthonormal extension $\tilde\V=\{v_1\ldots,v_{2m}\}$ of $\V$ such that the linear span of $\tilde\V$ contains $W$. By lemmas \ref{lem2.2}, \ref{lem2.3}, \ref{lem2.4} it follows that
\bes
[M(\V,r),M]\subseteq\tilde\V^{\perp}\subseteq W^{\perp}.
\ees
Since the symmetric bilinear form on $V$ is nondegenerate, we conclude that
\bes
[M(\V,r),M]=(0).
\ees
Choose arbitrary elements $u_1,\ldots,u_k\in V$. Let $\tilde\V$ be an orthonormal extension of $\V$ such that the linear span of $\tilde\V$ contains $u_1,\ldots,u_k$.  Then $M(\V,r)\subseteq \sum_{t=1}^{m-1}M(\tilde\V,t)$ and therefore
$$
M(\V,r)R(u_1)\cdots R(u_k)\subseteq \sum_{t=1}^{m-1}M(\tilde\V,t).
$$
Hence
$$
[M(\V,r)R(u_1)\cdots R(u_k),M]\subseteq \sum_{t=1}^{m-1}[M(\tilde\V,t),M]=0.
$$
We proved that the $A$-bimodule generated by $M(\V,r)$ is an ideal in $J$. Hence $M(\V,r)=(0)$.
This completes the proof of the lemma.

\ctd

{\bf\underline{Proof of Theorem 2.}}\\[2mm]
We proved that $M$ is a sum of regular and 1-dimensional bimodules over $F\cdot 1+\sum_{i=1}^{2n} Fv_i$.  Let
$$
M=M(1)+\sum_{i=1}^{2n}M(v_i)+M(v_1\cdots v_{2n}),
$$
where an arbitrary element from $M(1)$ is $\sim 1$ in $\V=\{v_1,\ldots,v_{2n}\}$, an arbitrary element from $M(v_i)$ is $\sim v_i$, an arbitrary element $x\in M(v_1\cdots v_{2n})$ is $\sim v_1\cdots v_{2n}$, hence $x\circ v_i=0,\,1\leq i\leq 2n$.
\smallskip

Clearly, $[M(1),M(1)]\subseteq F\cdot 1$.  We consider the skew-symmetric bilinear form
\bes
\g:M(1)\times M(1)\rightarrow F,\ \g(x\,|\,y)=[x,y].
\ees
For an arbitrary element $x(i)\in M(v_i)$ let $x'(1)=x(i)\circ v_i,\ x(i)=x'(1)\circ v_i$.  Hence 
$$
M(v_i)=M(1)\circ v_i,
$$
 and 
$$
[x,y\circ v_i]=\g(x\,|\,y)v_i  \hbox{ for } x,y\in M(1).
$$

Let $1\leq i\neq j\leq 2n;\ x,y\in M(1)$. Then $[x\circ v_i,y\circ v_j]\sim v_iv_j$,  hence 
$$
[x\circ v_i,y\circ v_j]=0.
$$
Similarly,
$$
[x\circ v_i,y\circ v_i]=[x,y\circ v_i]\circ v_i=\g(x\,|\,y).
$$
Let $W\subseteq M(1)$ be a finite-dimensional subspace of $M(1)$. Then
$$
J(W)=(F\cdot 1+\sum_{i=1}^{2n}Fv_i)+(W+\sum_{i=1}^{2n}W\circ v_i)
$$ 
is a subsuperalgebra of $J$.  Let 
$$
W^{\perp}=\{w\in W\,|\,\g(w\,|\,W)=(0)\}.
$$
Then $W^{\perp}+\sum_{i=1}^{2n} W^{\perp}\circ v_i$ is an ideal in $J(W)$ and $J(W)/(W^{\perp}+\sum_{i=1}^{2n} W^{\perp}\circ v_i)$ is a simple finite-dimensional Jordan superalgebra.  
It follows from  classification of simple finite-dimensional Jordan superalgebras  that $W=W^{\perp}$, $[M(1),M(1)]=(0)$.

This implies that 
\bes
[M,M]&=&[M(1),M(v_1\cdots v_{2n})]+[M(v_1\cdots v_{2n}),M(v_1\cdots v_{2n})]\\
&\subseteq& F\cdot 1+\{v_1,\ldots,v_{2n}\}^{\perp}.
\ees
Extending the orthonormal system $v_1,\ldots,v_{2n}$, we get
$$
[M,M]\subseteq F\cdot 1,
$$
and 
$$
[M(1)+\sum_{i=1}^{2n}M(v_i),M(1)+\sum_{i=1}^{2n}M(v_i)]=(0).
$$
The subspace $[M(1)+\sum_{i=1}^{2n}M(v_i),M(v_1\cdots v_{2n})]$ has zero intersection with $F\cdot 1$, hence
$$
[M(1)+\sum_{i=1}^{2n}M(v_i),M]=[M(1)+\sum_{i=1}^{2n}M(v_i),M(v_1\cdots v_{2n})]=(0).
$$
Let us show that for $n\geq 2$
$$
(F\cdot 1+\sum_{i=1}^{2n}Fv_i)D(M,M)=(0).
$$
We have
\bes
&(F\cdot 1+\sum_{i=1}^{2n}Fv_i)D(M(1)+\sum_{i=1}^{2n}M(v_i),M(1)+\sum_{i=1}^{2n}M(v_i))&\\
&\subseteq [M(1)+\sum_{i=1}^{2n}M(v_i),M(1)+\sum_{i=1}^{2n}M(v_i)]=(0).&
\ees
For $n\geq 2$
\bes
(F\cdot 1+\sum_{i=1}^{2n}Fv_i)D(M(1)+\sum_{i=1}^{2n}M(v_i),M(v_1\cdots v_{2n}))
\ees
has zero intersection with $F\cdot 1$,  hence 
\bes
(F\cdot 1+\sum_{i=1}^{2n}Fv_i)D(M(1)+\sum_{i=1}^{2n}M(v_i),M(v_1\cdots v_{2n}))=(0).
\ees
Finally, 
\bes
(F\cdot 1+\sum_{i=1}^{2n}Fv_i)D(M(v_1\cdots v_{2n}),M(v_1\cdots v_{2n}))=(0).
\ees
because $v_i\cdot M(v_1\cdots v_{2n})=(0)$. 
Extending the orthonormal system we get
$$
AD(M,M)=(0).
$$
If $x\in M$ and $[x,M]=(0)$ then for an arbitrary element $a\in A$ we have
$$
[x\circ a,M]\subseteq [x,a\circ M]+aD(x,M)=(0). 
$$
The subspace
$$
\{ x\in M\,|\, [x,M]=(0)\} 
$$ 
is an ideal of the superalgebra $J$. Since $[M,M]\neq (0)$ we conclude that this ideal is $(0)$ and therefore 
$$
M(1)+\sum_{i=1}^{2n}M(v_i)=(0),\ M=M(v_1\cdots v_{2n}),\ M\circ v_i=(0),\, 1\leq i\leq 2n.
$$
Extending the orthonormal system we get
$$
M\circ V=(0).
$$
This completes the proof of Theorem 2.

\ctd

\section{Proof of Theorem 3}

\hspace{\parindent}
Let $J=A+M$ be a simple Jordan superalgebra, $\dim_FA=\infty,\ A=J(V_1)\oplus J(V_2)$ is a direct sum of two simple Jordan algebras of symmetric nondegenerate bilinear forms.

Let $e,f$ be the identity elements of the Jordan algebras $J(V_1),\,J(V_2)$ respectively.
The decomposition
$$
M=\{e,M,e\}+\{e,M,f\}+\{f,M,f\}
$$
is the Pierce decomposition of $M$ with respect to the idempotents $e,f$.
\begin{lem}\label{lem3.1}
$M=\{e,M,f\}$.
\end{lem}
\prf
Let $P_e$ and $P_f$ be projections of the space $[\{e,M,f\},\{e,M,f\}]$ to $\{e,A,e\}$ and $\{f,A,f\}$ respectively.,
$$
[\{e,M,f\},\{e,M,f\}]\subseteq  P_e+P_f.
$$
Let us show that $P_f$ is an ideal in the algebra $\{f,A,f\}$.  
Observe that in a Jordan algebra $J$ with orthogonal idempotents $e,f$  whose sum is 1, the following inclusion holds
$$
\{\{e,J,f\},\{f,J,f\},\{e,J,f\}\}\subseteq \{e,J,e\}.
$$
In fact,  in a special algebra $J$  for any $a,b\in J$ we have
$$
\{\{e,a,f\},\{f,b,f\},\{e,a,f\}\}=\tfrac14\{e,\{a,\{f,b,f\},a\},e\}.
$$
Since $f=1-e$ and this equation is linear on $b$,  by the Shirshov-Macdonald theorem it holds in any Jordan algebra.
Linearizing it by $a$, we get the desired  inclusion.

Therefore,  in our superalgebra $J=A+M$ we have 
\bee\label{id3.1}
\{\{e,J,f\},\{f,J,f\},\{e,J,f\}\}\subseteq \{e,J,e\}.
\eee
In particular,  
\bes
\{\{e,M,f\},\{f,A,f\},\{e,M,f\}\}&\subseteq &\{e,A,e\}=J(V_1),\\
\{\{e,M,f\},\{f,M,f\},\{e,M,f\}\}&\subseteq &\{e,M,e\}.
\ees
Now for arbitrary elements $x,y\in \{e,M,f\},\, a\in \{f,A,f\},\, z\in \{f,M,f\}$ we have
\bes
[x,y]\cdot a&=&\{x,a,y\}-[x\cdot a,y]-[x,a\cdot y]\in \{e,A,e\}+P_f,\\ \ 
[x,y]\cdot z&=&-\{x,z,y\}+[x,z]\cdot y+x\cdot [z, y]\in \{e,M,e\},
\ees
which implies that $P_f$ is an ideal in $\{f,A,f\}=J(V_2)$ and $P_f\cdot \{f,M,f\}=(0)$.
Similarly,  $P_e$ is an ideal in $J(V_1)$ and $P_e\cdot \{e,M,e\}=(0)$.
Since the algebra $J(V_2)$ is simple,  if  $P_f\neq (0)$ then $f\in P_f$ and $\{f,M,f\}=(0)$.
Similarly,  if $P_e\neq (0)$ then $\{e,M,e\}=(0)$.  
Therefore,  if both $P_e\neq (0)$ and $P_f\neq (0)$ then $M=\{e,M,f\}$ and the lemma is proved.

Let $P_e=(0),\, P_f\neq (0)$.  Then $\{f,M,f\}=(0)$ and $[\{e,M,f\},\{e,M,f\}]\subseteq \{f,A,f\}$, which shows that 
$\{f,A,f\}+\{e,M,f\}$ is a nonzero ideal in $J$,  a contradiction.

Finally,  if  $[\{e,M,f\},\{e,M,f\}]=(0)$  then $\{e,M,f\}$ is an ideal of $J$, a contradiction again.
This completes the proof of the lemma.

\ctd

\begin{lem}\label{lem3.2}
Let $Fe+M$ be a Jordan superalgebra, $e^2=e,\, e\cdot x=\tfrac12 x$ for an arbitrary $x\in M$. Then $[M,M]=(0)$ or $\dim M=2$.
\end{lem}
\prf
Suppose that $x,y\in M,\,[x,y]\neq 0$. For any $z\in M$ we have by\eqref{Jordan_operator}
\bes
&z(R(x)R(e)R(y)-R(y)R(e)R(x)+R([x,y]\cdot e))&\\
&=z(R(x\cdot e)R(y)-R(y\cdot e)R(x)+R([x,y])R(e)).&
\ees
Hence
\bes 
[z,x]y-[z,y]x+[x,y]z=\tfrac12 [z,x]y-\tfrac12 [z,y]x+\tfrac12 z[x,y],
\ees
and therefore
\bes
[x,y]z=[z,y]x-[z,x]y.
\ees
The element $z$ is a linear combination of $x$ and $y$. This completes the proof of the lemma.

\ctd

\begin{lem}\label{lem3.3}
Let $\dim_FV\geq 2$ and let $J(V)$ be a simple Jordan algebra of a nondegenerate bilinear form on $V$, $e$ is the identity element of $J(V)$. Let $J=J(V)+M$ be a Jordan superalgebra, $e\cdot x=\tfrac12 x$ for an arbitrary element $x\in M$. Then $[M,M]=0$.
\end{lem}
\prf
Let $V\times V\rightarrow F,\ v\times w\mapsto (v\,|\,w)$ be a nondegenerate symmetric bilinear form that defined $J(V)$. Choose elements $v_1,v_2\in V$ such that $(v_1\,|\,v_1)=(v_2\,|\,v_2)=1,\ (v_1\,|\,v_2)=0$.
 Then $e_1=\tfrac12(e+v_1),\ e_2=\tfrac12(e-v_1)$ are orthogonal idempotents, $e_1+e_2=e$. We have
  $$
  M=M_1\oplus M_2,\ M_i=M\cdot e_i,\,i=1,2.
  $$
Consider the subsuperalgebra $Fe_1+M_1$. By lemma \ref{lem3.2}, $\dim_FM_1=2$ or $[M_1,M_1]=(0)$.  Suppose that $\dim_FM_1=2$. The element $v_2$ lies in $\{e_1,J(V),e_2\}$ and $v_2^2=e$. Hence $M_2=M_1\cdot v_2$. We showed that $\dim_FM_2=2,\, \dim_FM=4$.

There exists a finite-dimensional subspace $W\subseteq V$ such that $[M,M]\subseteq W$, the restriction of the form on $W$ is nondegenerate and $\dim_FW\geq 2$. If $[M,M]\neq (0)$ then the superalgebra $Fe+W+M$ has a nonunital simple homomorphic image with nonzero odd part and even part of dimension$\geq 3$. Such simple superalgebras do not exist (see \cite{Kac, RZ}).
\smallskip

Supppose that $[M_1,M_1]=(0)$.  Then 
\bes
[M,M]&=&[M_1,M_2]+[M_2,M_2]\subseteq Fe_2+\{e_1,A,e_2\},\\
\{e_1,A,e_2\}&=&v_1^{\perp}=\{v\in V\,|\, (v_1\,|\, v)=0\}.
\ees
The subspace $[M,M]$ is invariant with respect to inner derivations from $D(V,V)$. Let $[M,M]\ni a\neq 0,\, a=\a(e-v_1)+v,\, (v_1\,|\,v)=0,\,\a\in F$.

Suppose that $v\neq 0$. Then there exists an element $v'\in v_1^{\perp}$ such that $(v\,|\,v')=1$. Then $v_1D(v',v_1)=-v',\, vD(v',v_1)=v_1$.

Now, $aD(v',v_1)=\a v'+v_1=e_1-e_2+\a v'\not\in Fe_2+\{e_1,A,e_2\}$.

Suppose that $[M,M]\subseteq F(e-v_1)$. Choose $0\neq v'\in v_1^{\perp}$. Then 
$$
(e-v_1)D(v',v_1)=v'\not\in F(e-v_1).
$$
This complets the proof of the lemma.

\ctd

\begin{lem}\label{lem3.4}
Let $J=A+M$ be a Jordan superalgebra, $A=Fe_1\oplus Fe_2;\ e_1,e_2$ are orthogonal idempotents, $M=\{e_1,M,e_2\}$. Then there exists a subspace $M'\subseteq M$ of codimension $\leq 2$ such that $[M',M]=(0)$.
\end{lem}
\prf
1) Suppose at first that the superalgebra $J$ is finite dimensional and simple. Then from the classification (see \cite{Kac, RZ}) it follows that $J\cong D_t,\,\dim_FM=2$.
\smallskip

2) Next suppose that the superalgebra $J$ is finite-dimensional but not necessary simple.  Let $I\triangleleft J,\,I\neq J$ be a maximal ideal of $J$.

Suppose that $I\cap A\neq (0)$.  Then one of the idempotents $e_i$ lies in $I$, the other one does not. Let $e_1\in I$. This implies that $M=\{e_1,M,e_2\}\subseteq I$, and therefore $I=Fe_1+M,\, [M,M]\subseteq Fe_1$. By lemma \ref{lem3.2} $[M,M]=0$ or $\dim_FM=2$, which completes the proof in this case.

Let $I\cap A=(0)$. Since the superalgebra $J/I$ is simple, it follows from 1) that the codimension of $I$ in $M$ is $\leq 2$ and $[I,M]\subseteq I\cap A=(0)$.
\smallskip

3) Now let's drop the assumption that $\dim_FJ<\infty$. The superalgebra $J$ is locally finite-dimensional. Now the result follows from 2). This completes the proof of the lemma.

\ctd

\begin{lem}\label{lem3.5}
Let $J=A+M$ be a Jordan superalgebra, $A=Ff+J(V),\ f^2=f,\ J(V)$ is a Jordan algebra of a nondegenerate symmetric form on $V$, $\dim V=\infty, \ M=\{f,M,e\}$, where $e$ is the identity of $J(V)$. Then $[M,M]=(0)$.
\end{lem}
\prf
Let $v_0\in V,\, (v_0\,|\,v_0)=1,\ e_1=\tfrac12(e+v_0),\ e_2=\tfrac12 (e-v_0),$ then
$$
\{f+e_1,J,f+e_1\}=(Ff\oplus Fe_1)+M_1,\  M_1=\{f,M,e_1\}.
$$
By lemma \ref{lem3.4} there exists a subspace $M_1'\subseteq M_1$ of codimension $\leq 2$ such that $[M_1',M_1]=(0)$.
Define the skew-symmetric bilinear form $\langle ,\rangle: M\times M\rightarrow F$,  if $[x,y]=\a f+a,\, a\in J(V)$, then $\langle x\,|\,y\rangle=\a$.  Clearly, $\la M_1\,|\,M_2\ra=(0)$, hence $\la M_1'\,|\,M\ra=(0)$.  

Similarly we find a subspace $M_2'\subseteq M_2,\, [M_2',M_2]=(0),\ |M_2:M_2'|\leq 2$.
Let $M'=M_1'+M_2'$, then $\la M'\,|\,M\ra=(0),\ |M:M'|\leq 4$.
\smallskip

For arbitrary elements $x,y\in M$  we have $fD(x,y)=0$, hence 
$$
AD(x,y)\subseteq J(V).
$$
Therefore, for an arbitrary element $a\in J(V)$ we have
$$
[x\cdot a,y]+[y\cdot a,x]\in J(V),
$$
or, equivalently,
$$
\la x\cdot a\,|\,y\ra=\la x\,|\, y\cdot a\ra.
$$
Thus a multiplication by an element $a\in J(V)$ is a symmetric operator with respect to the form $\la x\,|\,y\ra$. This implies that
$$
\tilde M=\{x\in M\,|\,\la x\,|\,M\ra=(0)\}
$$
is an $A$-bimodule, $M'\subseteq \tilde M,\ |M:\tilde M|\leq 4,\ [\tilde M,\tilde M]\subseteq J(V)$.
\smallskip
Hence $J(V)+\tilde M$ is a subsuperalgebra of $J$, and by Lemma \ref{lem3.3} 
$$
[\tilde M,\tilde M]=(0).
$$
Choose  arbitrary elements $x_1\in M_1,\,y_2\in M_2$. The subspace 
$$
\tilde V=\{ v\in V\,|\, x_1\cdot v\in \tilde M, \, y_2\cdot v\in \tilde M,\ (v_0\,|\,v)=0\}
$$
has finite codimension in $V$. Since $\dim_FV=\infty$, the subspace above contains an element $v'$ such that $v'^2=e$,  it follows from nondegeneracy of the form.   Observe that $M_2\cdot v'\in \{f,M,e_2\}\cdot\{e_1,J,e_2\}\subseteq \{f,M,e_1\}=M_1$. Similarly,  $M_1\cdot v'\subseteq M_2$.

 It  follows from \eqref{SJ2} that 
$$
((y_2v')\cdot v')\cdot e_1+((e_1v')\cdot v')y_2=2(y_2v')(v'e_1),
$$
which implies $(y_2\cdot v')\cdot v'=\tfrac14y_2$.  

Now,  let $g,h$ be orthogonal idempotents in a Jordan algebra $J$, then we have $JU(g,h)U(h)=(0)$ and hence 
$$
JU(g,h)U(J(U(h)))\subseteq JU(g,h)U(h)U(J)U(h)=(0)
$$
 or $\{\{h,J,h\},\{g,J,h\},\{h,J,h\}\}=(0)$.
Taking $h=f+e_1,\,g=e_2$ in the superalgebra $J$, we have  $M_1\subseteq \{h,J,h\},\ v'\in \{g,J,h\}$, hence $\{M_1,v',M_1\}=(0)$
and 
$$
0=\{x_1,v',y_2\cdot v'\}=[x_1\cdot v',y_2\cdot  v']-[(y_2\cdot v')\cdot v',x_1]-[x_1,y_2\cdot v']\cdot v''.
$$
We have  $[x_1,y_2\cdot v']\in[M_1,M_1]\subseteq Ff+Fe_1, \ (Ff+Fe_1)\cdot v'\subseteq Fv'$.  Also, $[x_1\cdot v',y_2\cdot  v']\in [\tilde M,\tilde M]=(0)$, which implies 
$$
[x_1,y_2]\in Fv'.
$$
There is another element $v''\in\tilde V$, such that $(v'')^2=e,\ v''$ is not collinear with $v'$.  Arguing as above,  we get  $[x_1,y_2]\in Fv''$. Hence $[x_1,y_2]=0,\ [M_1,M_2]=(0)$.
\smallskip

Now we have
$$
[M,M]=[M_1,M_1]+[M_2,M_2]\subseteq Ff+Fe+Fv_0.
$$
Choosing another element $v_0'$ with the same properties, we get
$$
[M,M]\subseteq Ff+Fe+Fv_0'.   
$$
Therefore, 
$$
[M,M]\subseteq Ff+Fe.
$$
Let $x,y\in M,\ [x,y]=\a f+\b e,\, \b\neq 0$.  The subspace
$$
\{v\in V\,|\, [x\cdot v,y]=[x,y\cdot v]=0\}
$$
has a finite codimension in $V$. Let $v$ be a nonzero element from this subspace. Then by \eqref{id3.1}
$$
\{x,v,y\}=[x\cdot v,y]-[y\cdot v,x]-[x,y]v\in \{f,J,f\}=Ff.
$$
This implies $[x,y]\cdot v=\b v=0$,  a contradiction.  Hence $[M,M]\subseteq Ff$.

We showed that $\tilde M=\{x\in M\,|\,\la x,M\ra=(0)\}$ is an $A$-subbimodule of $M$.  Since $[\tilde M,M]=\la \tilde M,M\ra=(0)$,  it follows that $\tilde M$ is an ideal of $J$.   Suppose that $\tilde M\neq M$ and factor the ideal $\tilde M$ out.  Then without loss of generality we can assume that $M\neq (0)$ and $\dim_FM\leq 4$. 

The subspace $\{v\in V\,|\,M\cdot v=(0)\}$ has finite codimension in $V$, hence it contains an element $v$ such that $(v,v)=1$.  For an arbitrary element $x\in M$ we have  
$$
4(x\cdot v)\cdot v=x=0.
$$
We showed that $\tilde M=M,\ [M,M]=(0)$.  This completes the proof of the lemma.

\ctd

{\bf\underline{Proof of Theorem 3.}}\\[2mm]
Suppose that $A=J(V_1)\oplus J(V_2); \ \dim_F V_i\geq 2;\, i=1,2.$
The algebra $A$ is infinite dimensional,  hence at least one of the spaces $V_1,V_2$ is infinite dimensional. 
Suppose that $\dim_FV_2=\infty$.  Let $f,e$ be the identity elements of the algebras $J(V_1),J(V_2)$, respectively.  
By Lemma \ref{lem3.1}, $M=\{e,M,f\}$.  
The identity $f$ is a sum of two orthogonal idempotents, $f=f_1+f_2$.  The subsuperalgebra $J_1=\{f_1+e,J,f_1+e\}$ satisfies the assumptions of Lemma \ref{lem3.5},  $(J_1)_{\0}=Ff_1+J(V_2),\ (J_1)_{\1}=M_1=\{ f_1,M,e\}$. By Lemma \ref{lem3.5} $[M_1,M_1]=(0)$.
Consider also the subsuperalgebra $J_2=\{f_2+e,J,f_2+e\}=Ff_2+J(V_2)+M_2,\ M_2=\{f_2,M,e\},\ M=M_1+M_2$.   As before, we have $[M_2,M_2]=(0)$,  hence $[M,M]=[M_1,M_2]=\{f_1,J(V_1),f_2\}\subseteq J(V_1)$.

The superalgebra $J(V_1)+M$  satisfies the assumptions of Lemma \ref{lem3.3}, hence by  this lemma $[M,M]=(0)$,  which contradicts the simplicity of the superalgebra $J$.   This completes the proof.

\ctd

\bigskip

The authors thank the referee for various usefull remarks.

 \end{document}